%
%
\magnification=\magstephalf
\input amstex
\documentstyle{amsppt}
\NoBlackBoxes
\TagsOnRight
\loadbold


\def\qi{\bold i}
\def\qj{\bold j}
\def\qk{\bold k}
\def\qxi{\boldsymbol \xi}
\def\qpi{\boldsymbol \pi}
\def\qtau{\boldsymbol \tau}
\def\qalpha{\boldsymbol \alpha}
\def\qlambda{\boldsymbol \lambda}
\def\qchi{\boldsymbol \chi}
\def\qfone{\bold f^{(1)}}
\def\qalfone{\boldsymbol \alpha^{(1)}}
\def\locqed{\qed\medskip}          

\settabs  12 \columns

\topmatter
\title A Quaternionic Proof of the Representation Formula of a Quaternary 
       Quadratic Form\endtitle
\rightheadtext{Quaternionic Proof of Quadratic Form Representation Formula}

\author Jesse Ira Deutsch \endauthor

\address Mathematics Department, University of Botswana, Private Bag 0022,
Gaborone, Botswana \endaddress


\email deutschj\@mopipi.ub.bw\endemail


\subjclass Primary 11D09, 11D57, 11P05 \endsubjclass

\abstract 
The celebrated Four Squares Theorem of Lagrange states that 
every positive integer is the sum of four squares of integers.
Interest in this Theorem has motivated a number of different
demonstrations.
While some of these demonstrations prove the existence of 
representations of an integer as a sum of four squares, 
others also produce the number of such representations.
In one of these demonstrations, Hurwitz was able to use a
quaternion order to obtain the formula for the number
of representations.
Recently the author has been able to use certain quaternion
orders to demonstrate the universality of other quaternary
quadratic forms besides the sum of four squares.
In this paper we develop results analogous to Hurwitz's 
above mentioned work by delving into the number theory of 
one of these quaternion orders,
and discover an alternate proof of the representation
formula for the corresponding quadratic form.
\endabstract




\keywords Quaternions, quadratic form, representation formulas\endkeywords

\endtopmatter

\document

\head 1. Introduction\endhead

Four variable quadratic forms have been of interest since Lagrange
proved the theorem that every positive integer is the sum
of four integer squares.
A quadratic form over $\Bbb Z$ is called universal if
it represents each positive integer at least once.
It has been shown that there are a total of 54
quaternary quadratic forms without cross product terms
that are universal (see Duke [2] and Ramanujan [8]).
Hurwitz was able to demonstrate the universal property for the
quadratic form that is the sum of four squares by use of
a special quaternion order.
Hurwitz's study of the number theory of this order additionally
yielded the formula for the number of representations of any
positive integer as a sum of four squares.
In previous work, the author has been able to show the
universal property of seven other quadratic forms by consideration
of appropriate quaternion orders.
Here we develop the number theory of the quaternion order
related to $x^2 + y^2 + 2\, z^2 + 2\, w^2$, thus creating an
alternate proof of the representation formula and the
universality for this quadratic form.
In addition we find formulas for the number of representations
of four and eight times an odd number by this quadratic form
with certain restrictions on the parity of the variables
$x$, $y$, $z$, and $w$.
See Hurwitz [6] and Deutsch [1] for further background.

Using the notation of Deutsch [1], we let
$\qi,\, \qj,\, \qk$ be the standard noncommutative basis elements
of the quaternions.
Bold type is used to represent quaternions, so a typical one
is of the form
$$
\bold q  
     \ = \ 
  q_1 + q_2 \qi + q_3 \qj + q_4 \qk \,, \qquad
    q_1, q_2, q_3, q_4 \,\in {\Bbb R} 
\tag 1.1
$$
To each such $\bold q$ there corresponds a conjugate quaternion
denoted $\bold{\overline q}$ which has the value
$q_1 - q_2 \qi - q_3 \qj - q_4 \qk$.
The quaternionic norm, $N( \bold q )$ is simply the product 
$\bold q \cdot \bold {\overline q}$.

Recall that a ring of quaternions is
norm Euclidean if there exists a nonnegative $\delta < 1$ such that
for any quaternion $\bold q$
there exists a quaternion $\bold a$ in the ring for which
$N (\bold q - \bold a) \leq \delta$.

In general, a unit of a ring with unity is an invertible
element of the ring.
An associate of an element of the ring is that element multiplied
by a unit on the left or right.

For norm Euclidean rings of quaternions, there exist left and
right greatest common divisors of any two nonzero elements.
These are unique up to associates.
The right greatest common divisor can be written as a right linear
combination of the two nonzero elements under consideration.
An analogous result holds for left greatest common divisors.
For details and demonstrations, refer to Hardy and Wright [5, Ch.~XX].

The relevant quaternion order for the quadratic form
$x^2 + y^2 + 2 z^2 + 2w^2$ is denoted $H_{1,2,2}$ .
It is the $\Bbb Z$ module with generators
$$
\bold v_1 \ = \ 1,\ 
  \bold v_2 \ = \ \qi,\ 
  \bold v_3 \ = \ {1 \over 2} \left( 1 + \qi + \sqrt 2 \, \qj \right),\ 
  \bold v_4 \ = \ {1 \over 2} \left( 1 + \qi + \sqrt 2 \, \qk \right)
  \,.
\tag 1.2
$$
This module is a norm Euclidean ring with a total of 24 units,
namely
$$
\pm \left\{
\eqalign {
  &\bold v_1,\, \bold v_2,\, 
   \bold v_3,\, \bold v_4,\,
\cr
  &\bold v_3 - \bold v_1,\, \bold v_3 - \bold v_2,\,
   \bold v_4 - \bold v_3,\,
   \bold v_4 - \bold v_1,\,
   \bold v_4 - \bold v_2,\,
\cr
    &\bold v_3 - \bold v_2 - \bold v_1,\,
     \bold v_4 - \bold v_2 - \bold v_1,\,
     \bold v_4 + \bold v_3 - \bold v_2 - \bold v_1
\cr
}
\right\} 
\tag 1.3
$$
(see Deutsch [1]).
$H_{1,2,2}$ is closed under conjugation, and the norm maps $H_{1,2,2}$
into the nonnegative rational integers.
In addition to $H_{1,2,2}$, we shall also make use of the $\Bbb Z$
module generated by $\{ 1,\, \qi,\, \sqrt 2 \,\qj,\, \sqrt 2 \,\qk \}$.
This module will be denoted $H_{1,2,2}^0$.

For related notation and further properties of the above quaternion 
order see Deutsch [1].
More background material on quaternion orders and algebras is available 
in Pierce [7], Scharlau [9] and Vign\'eras [10].
An alternate proof of the representation formula using analytic
techniques is given in Fine [3].

\bigskip
\head 2. Residues modulo powers of $1 + \qi$ \endhead

In Hurwitz's study of the quaternion order for the sum of
four squares, elements of norm $2$ were very significant.
In particular, $1 + \qi$ had the special property that every
element of the order with norm divisible by $2$ was also
divisible by $1 + \qi$ on both the left and the right.
This property carries over to the case of $H_{1,2,2}$ with
some work.
Of course $1 + \qi = \bold v_1 + \bold v_2 \in H_{1,2,2}$.
We start by noting that 
$N ( \bold v_1 + \bold v_3 - \bold v_4 ) =  2$.
Also, in $H_{1,2,2}$ we have
$$
\eqalign{
\bold v_1 + \bold v_3 - \bold v_4 
  \ &= \ 
    \left( 
         {1 \over 2} - {1 \over 2} \qi + {{\sqrt 2} \over 2} \qj
    \right)
    \cdot (1 + \qi)
  \cr
  \ &= \
     (1 + \qi) \cdot
     \left( 
        {1 \over 2} - {1 \over 2} \qi - {{\sqrt 2} \over 2} \qk
     \right)
     .
  \cr
}
\tag 2.1
$$

\proclaim{Lemma 1}
The coset representatives of $H_{1,2,2}$ modulo $(1 + \qi)$ are
$\{ 0,\, 1,\, \bold v_3,\, 1 + \bold v_3 \}$.
This holds whether $(1 + \qi)$ is considered a left or
a right ideal.
\endproclaim

\demo{Proof}
Consider any 
$\bold g = g_1 \bold v_1 + g_2 \bold v_2 +g_3 \bold v_3 + g_4 \bold v_4$
in $H_{1,2,2}$.
Using equation (2.1), and letting $\bold h$ be an appropriate element in 
this order 
that may be different on each line, we have 
$$
\eqalign{
    g \ &= \ g_1 + g_2 \,\qi + g_3 \,\bold v_3 + g_4 \,(\bold v_3 + 1)
           + (1 + \qi) \,\bold h
    \cr	   
       &= \ (g_1 + g_4) + g_2 \,\qi + \bold v_3 \,(g_3 + g_4) 
           + (1 + \qi) \,\bold h
    \cr
       &= \ (g_1 + g_2 + g_4) + \bold v_3 \,(g_3 + g_4)
           + (1 + \qi) \,\bold h
	   \,.
    \cr
}
\tag 2.2
$$
The last line follows from the elementary identity
$\qi - 1 = \qi (1 + \qi) = (1 + \qi) \qi$.
Since $1 + \qi$ divides $2$, the only possibilities for the
integer quantities $g_1 + g_2 + g_4$ and $g_3 + g_4$ are
$0$ and $1$.
These possibilities generate the coset representatives of
the Lemma.
These coset representatives are not divisible by $1 + \qi$
since the norm of the differences of any two of them is
odd, while the norm of $1 + \qi$ is even.

The same argument hold in the case where $\bold h$ is on
the left in (2.2).
\locqed
\enddemo

\proclaim{Lemma 2}
The coset representatives of $H_{1,2,2}$ modulo $(1 + \qi)$ can
be chosen as
$\{ 0,\, 1,\, \bold v_3,\, \bold v_3^2 \}$.
\endproclaim

\demo{Proof}
This follows from
$$
\bold v_3^2 \ = \ 
  \bold v_3 - 1 \ = \ 
  \bold v_3 + 1 - 2 \ = \ 
  \bold v_3 + 1 - (1 + \qi) \cdot (1 - \qi) \ = \ 
  \bold v_3 + 1 - (1 - \qi) \cdot (1 + \qi)
\,.
\quad\qed
\tag 2.3
$$
\enddemo

We have already used the fact that if $1 + \qi$ divides an element
of $H_{1,2,2}$ then $2$ must divide the norm of that element.
It is nontrivial that the converse is also true.

\proclaim{Theorem 3}
Suppose $\bold g \in H_{1,2,2}$ with $N(\bold g )$ even.
Then $\bold g$ has $1 + \qi$ as a left factor in $H_{1,2,2}$.
Also $\bold g$ has $1 + \qi$ as a right factor in $H_{1,2,2}$.
\endproclaim

\demo{Proof}
Let us show that $\bold g$ has $1 + \qi$ as a right factor.
Considering the coset representatives of $H_{1,2,2}$ modulo
$1 + \qi$, $\bold g$ must have one of the forms
$$
0 + \bold h \, (1 + \qi) \,, \quad
1 + \bold h \, (1 + \qi) \,, \quad
\bold v_3 + \bold h \, (1 + \qi) \,, \quad
1 + \bold v_3 + \bold h \, (1 + \qi) \,.
\tag 2.4
$$
for some $\bold h \in H_{1,2,2}$.
Write $\bold g = \bold m + \bold h \, (1 + \qi)$.
Then
$$
\eqalign{
\bold g \cdot \bold{\overline g} 
  \ &= \  
    \left( \bold m + \bold h \,(1 + \qi) \right) \cdot
       \left( \bold {\overline m} + (1 - \qi) \,\bold{\overline h} \right)
  \cr
  \ &= \ 
    \bold m \cdot \bold{\overline m} 
      + \bold h \,(1 + \qi) \,\bold{\overline m}
      + \bold m \,(1 - \qi) \,\bold{\overline h}
      + 2 \, \bold h \cdot \bold{\overline h}
   \,.   
 \cr
}       
\tag 2.5
$$
For the center two terms, note that
$$
\overline{ \left( \bold h \, (1 + \qi) \, \bold{\overline m} \right)}
  \ = \  
  \bold m \, (1 - \qi) \, \bold{\overline h}
  \,,
\tag 2.6
$$
and that the sum of a quaternion and its conjugate is twice its
real coefficient when written in the standard basis, i.e.
for $\bold q$ as in (1.1) we have 
$\Re (\bold q) = \bold q + \bold{\overline q} = 2\, q_1$.
We must therefore find the real coefficient in the standard basis of
$\bold h \, (1 + \qi) \, \bold{\overline m}$.
Write $\bold h$ and $\bold m$ in the standard quaternion basis
$$
\bold h  
     \ = \ 
  h_1 + h_2 \qi + h_3 \qj + h_4 \qk \,, \qquad
\bold m  
     \ = \ 
  m_1 + m_2 \qi + m_3 \qj + m_4 \qk \,.
\tag 2.7
$$
Computation shows that
$$
\Re \left(
          \bold h \, (1 + \qi) \, \bold{\overline m}
    \right)
  \ = \ 
     (h_1 - h_2) m_1 + (h_1 + h_2) m_2 + (h_3+h_4) m_3
       + (h_4- h_3) m_4
  \,.        
\tag 2.8
$$
However $\bold m$ is restricted to the four values of the
coset representatives listed in Lemma 1.
We wish to show this real part must always be a rational integer.

If $\bold m = 0$ then     
$\Re \left(
          \bold h \, (1 + \qi) \, \bold{\overline m}
    \right) = 0$.
If $\bold m = 1$ then
$\Re \left(
          \bold h \, (1 + \qi) \, \bold{\overline m}
    \right) = h_1 - h_2$.
Write $\bold h$ in terms of the basis for $H_{1,2,2}$.
$$
\eqalign{
\bold h  
  \ &= \ 
     t_1 \bold v_1 + t_2 \bold v_2 + t_3 \bold v_3 + t_4 \bold v_4 
  \cr
  \ &= \ 
    \left( t_1 + {1 \over 2} t_3 + {1 \over 2} t_4 \right)
    + \left( t_2 + {1 \over 2} t_3 + {1 \over 2} t_4 \right) \qi
    + {{t_3} \over 2} \,\sqrt 2 \,\qj
    + {{t_4} \over 2} \,\sqrt 2 \,\qk
    \,.
  \cr  
}       
\tag 2.9
$$
A simple calculation shows that
$h_1 - h_2 = t_1 - t_2$ which is an element of $\Bbb Z$.

For $\bold m = \bold v_3$ we similarly find
$\Re \left(
          \bold h \, (1 + \qi) \, \bold{\overline m}
    \right) = t_1 + t_3 + t_4 \in \Bbb Z$.
In the case of $\bold m = 1 + \bold v_3$,
$\Re \left(
          \bold h \, (1 + \qi) \, \bold{\overline m}
    \right) = 2 t_1 - t_2 + t_3 + t_4 \in \Bbb Z$.

In all cases we find
$\Re \left(
          \bold h \, (1 + \qi) \, \bold{\overline m}
    \right)$
is a rational integer, so twice this quantity is an even integer.
Thus 
$$
N ( \bold g ) \ = \ \bold g \cdot \bold{\overline g}
  \ \equiv \ 
  \bold m \cdot \bold {\overline m}   \pmod 2
\tag 2.10
$$
in the rational integers.
As $\bold m$ varies through the cases listed above, the norm of
$\bold m$ takes on the values $0,\, 1,\, 1,\, 3$ in order.
Thus if $N(\bold g)$ is even, $\bold m$ must equal $0$.
Hence $\bold g = \bold h \, (1 + \qi)$ so we find $1 + \qi$
is a right factor of $\bold g$.
      
The demonstration that $1 + \qi$ is also a left factor of $\bold g$
is entirely analogous to the above proof.
\locqed
\enddemo

A simple induction argument yields the following.

\proclaim{Corollary 4}
Suppose $\bold g \in H_{1,2,2}$ and $2^s$ divides $N(\bold g )$ 
for some positive integer $s$.
Then $(1 + \qi)^s$ divides $\bold g$ on each of the left side
and the right side.
\endproclaim

Since every real number is in the centralizer of the quaternions,
we do not have to distinguish left and right cosets when considering
$H_{1,2,2}$ modulo a rational integer.

\proclaim{Lemma 5}
The set of coset representatives for $H_{1,2,2} \,/\, (2)$ can be
chosen as the twelve units listed in {\rm (1.3)} taken with the positive
sign, and the four nonunits 
$0,\, 1 + \qi,\, 1 + \bold v_3 + \bold v_4,\, 
   \qi + \bold v_3 + \bold v_4$.
\endproclaim

\demo{Proof}
A typical $\bold g \in H_{1,2,2}$ can be written
$\bold g = g_1 + g_2 \,\qi + g_3 \bold v_3 + g_4 \,\bold v_4$ 
with $g_1, \,\dots,\, g_4 \in \Bbb Z$.
Taken modulo $2$ we need only consider the values
$g_1,\, \dots,\, g_4 \in \{0,\, 1 \}$.
Shifting appropriate $g$'s by $2$ produces the listing of
the Lemma.
Note that $N(0) = 0$, $N(1 + \qi) = 2$, 
$N(1 + \bold v_3 + \bold v_4 ) = 6$,
$N(\qi + \bold v_3 + \bold v_4) = 6$.
Hence these elements cannot be units.

Since as free abelian groups 
$H_{1,2,2} \,/\, (2) \,\simeq\, 
  {\Bbb Z}^4 \,/\, 2 {\Bbb Z}^4 \,\simeq\,
  \left( {\Bbb Z} \,/\, 2 {\Bbb Z} \right)^4$,
a group of order $16$,
the Lemma lists the correct number of elements.
\locqed
\enddemo

\definition{Definition 6}
An element of $H_{1,2,2}$ is an odd quaternion if its norm is odd.
\enddefinition

\proclaim{Lemma 7}
Any odd quaternion in $H_{1,2,2}$ is congruent to a unit
modulo $2$.
\endproclaim

\demo{Proof}
Let $\bold b$ be odd.
If $\bold b$ were congruent to one of the non-units then we could
write
$\bold b = \bold r + 2 \bold h$ with $N(\bold r )$ even.
Thus $1 + \qi$ divides $\bold r$ on the left.
But $1 + \qi$ divides $2$ also, so $1 + \qi$ divides $\bold b$.
Hence the norm of $\bold b$ would be even, a contradiction.
\locqed
\enddemo

\proclaim{Lemma 8}
Multiplication by an odd quaternion on the left permutes the
twelve units of Lemma 5 modulo $2$.
The same holds for multiplication on the right.
\endproclaim

\demo{Proof}
Let $\bold b$ be odd, and let $\bold u$ be one of the above
mentioned twelve units modulo $2$.
Then $N(\bold b \, \bold u) = N(\bold b ) \, N(\bold u )$
is odd so $\bold b \, \bold u$  is congruent to one of these
twelve units modulo $2$.

Let $\bold u,\, \bold u^\prime$ be units and suppose
$\bold b \, \bold u \equiv \bold b \, \bold u^\prime \pmod 2$
in $H_{1,2,2}$.
Then $2 \,\big|\, \bold b (\bold u - \bold u^\prime)$.
Hence $4 \,\big|\, N(\bold b) \, N(\bold u - \bold u^\prime)$.
Thus $4 \,\big|\, N(\bold u - \bold u^\prime)$ as $\bold b$ is odd.
But this implies that
$(1 + \qi)^2 \,\big|\, \bold u - \bold u^\prime$.
Thus it follows that $2 \,\big|\, \bold u - \bold u^\prime$, so
$\bold u$ and $\bold u^\prime$ are equivalent modulo $2$.

Hence multiplication by odd $\bold b$ on the left induces a 
one-to-one onto mapping of the twelve units modulo $2$ to
themselves.
The same argument holds for multiplication on the right.
\locqed
\enddemo

\proclaim{Corollary 9}
Let $\bold b$ be an odd quaternion in $H_{1,2,2}$.
There exist two units $\bold u$ and $\bold u_1$
which satisfy the congruences
$$
\bold b \, \bold u
  \ \equiv \ 
    \bold u_1 \, \bold b
  \ \equiv \ 
    1 \pmod 2
  \,.    
\tag 2.11  
$$  
\endproclaim

\demo{Proof}
Since multiplication by $\bold b$ induces a permutation among
the units modulo $2$, and $1$ is a unit, the Corollary follows.
\locqed
\enddemo

\bigskip
\head 3. Primary quaternions \endhead

Hurwitz used the concept of primary quaternions to enable 
the counting of quaternions with requisite properties for 
the order related to the quadratic form
$x^2 + y^2 +z^2 + w^2$.
It turns out that the exact analogue of Hurwitz's definition
accomplishes the same goal for the order $H_{1,2,2}$.

\proclaim{Lemma 10}
Let $a \in \Bbb Z$, $a \geq 0$.
Then the left ideal generated by $2^a (1 + \qi)$ in 
$H_{1,2,2}$ is the same as the right ideal generated
by this element.
This set is the same as the two sided ideal generated
by $2^a (1 + \qi)$.
\endproclaim

\demo{Proof}
Let $\bold q = \bold h \cdot 2^a (1 + \qi)$.
Then $2^{2a + 1} \,\big|\, N(\bold q)$.
By a Lemma above, $(1 + \qi)^{2a+1}$ divides $\bold q$ on
the left.
Since $(1 + \qi)^2$ is $2 \qi$, we find there exists
$\bold h^\prime$ such that
$\bold q = 2^a (1 + \qi) \cdot \bold h^\prime$.
A similar argument holds in the other direction.

Any element in a one sided ideal is contained in the two
sided ideal.
In the situation of this Lemma, we again have that
$2^{2a + 1}$ divides the norm of any element of the two
sided ideal.
Thus $2^a (1 + \qi)$ divides any element of the two sided
ideal from the left or right.
\locqed
\enddemo

In particular, it makes sense to use the notation 
$\left( 2 (1 + \qi) \right)$ for any of the one sided 
or two sided ideals since they are all the same.

\proclaim{Lemma 11}
The set $\{ 1,\, 1 + 2 \, \bold v_3\}$ forms a multiplicative
subgroup of $H_{1,2,2}$ modulo the ideal 
$\left( 2 (1 + \qi) \right)$.
\endproclaim

\demo{Proof}
The only nontrivial product to check is
$(1 + 2 \bold v_3 )^2$.
Consider the difference
$(1 + 2 \bold v_3)^2 - 1$.
This equals $4 \cdot (\bold v_3 + \bold v_3^2 )$.
Since $2 (1 + \qi)$ divides $4$,
$(1 + 2 \bold v_3)^2$ is congruent to $1$ modulo
$2 (1 + \qi)$.
\locqed
\enddemo

\definition{Definition 12}
A primary quaternion is an element of $H_{1,2,2}$ which is
congruent to $1$ or $1 + 2 \bold v_3$ modulo the ideal
$\left( 2 (1 + \qi) \right)$.
\enddefinition

By Lemma 11 we immediately see that the product of two primary
quaternions is primary.

\proclaim{Lemma 13}
Primary quaternions are elements of $H_{1,2,2}^0$.
\endproclaim

\demo{Proof}
Let $\bold q$ be a primary quaternion.
Then $\bold q$ is congruent to $1$ or $1 + 2 \bold v_3$ modulo
$2 (1 + \qi)$.
In particular $\bold q$ can be written in the form
$1 + 2 \bold g$ for some $\bold g \in H_{1,2,2}$.
Thus
$$
\eqalign{
\bold q \ &= \ 
     1 + 2 (g_1 \bold v_1 + g_2 \bold v_2 + g_3 \bold v_3 
            + g_4 \bold v_4)
  \cr
  \ &= \ 
     1 + 2g_1 + 2g_2 \qi + g_3 (1 + \qi + \sqrt 2 \qj)
       + g_4 (1 + \qi + \sqrt 2 \qk)
  \cr
}    
\tag 3.1
$$
an element of $H_{1,2,2}^0$.
\locqed
\enddemo

\proclaim{Lemma 14}
The set of 24 right sided associates of an odd quaternion always
contains one which is primary.
The same holds for the left sided associates.
\endproclaim

\demo{Proof}
Start with an arbitrary quaternion congruent to $1$ modulo $2$.
It must be of the form $1 + 2 \,\bold g$, $\bold g \in H_{1,2,2}$.
To consider this quantity modulo $2 (1 + \qi)$ we need only
look at the possible values of $\bold g$ modulo $1 + \qi$, namely
$\{ 0,\, 1,\, \bold v_3,\, \bold v_3^2 \}$.
Thus the possible values of $1 + 2 \bold g$ are
$$
\eqalign{
1 + 2 (1 + \qi) \bold g^\prime \,, \quad
1 + 2  &\left( (1 + \qi) \bold g^\prime + 1 \right) \,, \quad
1 + 2  \left( (1 + \qi) \bold g^\prime + \bold v_3 \right) \,, \quad
  \cr
&1 + 2  \left( (1 + \qi) \bold g^\prime + \bold v_3^2 \right) 
  \cr
}
\tag 3.2
$$
for some $\bold g^\prime \in H_{1,2,2}$.
Taking this modulo $2 (1 + \qi)$ the possible values become
$1,\, 1 + 2 ,\,  1 + 2 \bold v_3 ,\, 1 + 2 \bold v_3^2$.
Since $2 (1 + \qi)$ divides $4$, we observe 
$1 + 2 = 3 \equiv -1 \pmod 4$ so
$3 \equiv -1 \pmod {2 (1 + \qi)}$.
Also
$$
1 + 2 \bold v_3^2 \ \equiv \ 
  1 + 2 (\bold v_3 - 1) \ \equiv \ 
  -1 + 2 \bold v_3 \ \equiv \ 
  -1 - 2 \bold v_3 \pmod 4
\,.
\tag 3.3
$$
Clearly these equivalences also hold modulo $2 (1 + \qi)$.
The possible values for $1 + 2 \bold g$ modulo
$2 (1 + \qi)$ reduce to 
$1,\, -1,\, 1 + 2 \bold v_3 ,\, -1 - 2 \bold v_3$.

Fix an odd quaternion $\bold b$.
There exists a unit $\bold u$ such that 
$\bold b \, \bold u \equiv 1 \pmod 2$.
Also $\bold b \, (-\bold u) \equiv -1 \equiv 1 \pmod 2$.
We conclude that one of $\bold b \, \bold u$ and
$\bold b \, (-\bold u)$ is congruent to a member of the
set $\{1 ,\, 1 + 2 \bold v_3 \}$ modulo $2 (1 + \qi)$.

For left sided associates, we need only consider
$\bold u_1 \, \bold b$ and $-\bold u_1 \, \bold b$.
\locqed
\enddemo

\proclaim{Lemma 15}
Exactly one of the right sided associates of an odd quaternion
is primary.
The same holds for the left sided associates.
\endproclaim

\demo{Proof}
Let $\bold b$ be as above.
Suppose $\bold b \, \bold u$
and $\bold b \, \bold u_1$ are both primary for $\bold u$
and $\bold u_1$ units.
Then $\bold b \, \bold u$ and $\bold b \, \bold u_1$ are
both congruent to $1$ modulo $2$ in $H_{1,2,2}$.
Let $\bold u_2 = \bold u^{-1} \, \bold u_1$.
Then $2$ divides $\bold b \, \bold u \, (1 - \bold u_2)$.
Taking norms, we see that in the rational integers $4$ must 
divide $N(\bold b) N(\bold u) N(1 - \bold u_2)$.
Thus $4$ evenly divides $N(1 - \bold u_2)$.
But $\bold u_2$ is a unit, and computation shows that of the
24 possible units, $4$ divides $N(1 - \bold u_2)$ only in the
case $\bold u_2 = \pm 1$.

If $\bold u_2 = 1$ then $\bold u = \bold u_1$ and the Lemma is
proven.
If $\bold u_2 = -1$ then $- \bold u = \bold u_1$ thus
$$
\bold b \, \bold u_1 \ \equiv \ 
  -1 \hbox{ or }  -1 -2 \bold v_3 \pmod {2 (1 + \qi) }
  \,.
\tag 3.4
$$
Computation shows that 
$\{1,\, 1 + 2 \bold v_3,\, -1,\, -1-2\bold v_3 \}$
are all distinct modulo $2 (1 + \qi)$ as the norms of all
pairs of differences are divisible by $4$ but never by $8$.
Thus $\bold b \, \bold u_1$ in the above equation can never
be primary.  

A very similar proof holds for the products $\bold u \, \bold b$ where
$\bold u$ is a unit.
\locqed
\enddemo

\proclaim{Lemma 16}
Let $\bold b$ be a primary quaternion.
If $\bold b \equiv 1 \pmod {2 (1 + \qi)}$
then the conjugate of $\bold b$, $\bold {\overline b}$
is also primary.
If $\bold b \equiv 1 + 2 \bold v_3 \pmod {2 (1 + \qi)}$
then $-\bold {\overline b}$ is also primary.
\endproclaim

\demo{Proof}
We start from the congruence of the left ideal
$$
\bold b \ \equiv \ 
  1 \hbox{ or } 1 + 2 \bold v_3
  \pmod{2 (1 + \qi)}
\tag 3.5
$$
and by taking conjugates we have the congruence of a right ideal
$$
\bold {\overline b} \ \equiv \ 
  1 \hbox{ or } 1 + 2 \bold {\overline v}_3
    \pmod {2 ( 1 - \qi)}
\,.
\tag 3.6
$$    
However $1 - \qi = -\qi ( 1 + \qi) = (1 + \qi) (-\qi)$
so we have
$$
\bold {\overline b} \ \equiv \ 
  1 \hbox{ or } 1 + 2 \bold {\overline v}_3
    \pmod {2 ( 1 + \qi)}
\,.
\tag 3.7
$$    
Thus we are back again in the case where the left, right and
two sided ideals are the same.
At this point we need only show
$$
1 + 2 \bold v_3 \ \equiv \ 
  - \left( 1 + 2 \bold {\overline v}_3 \right)
  \pmod {2 ( 1 + \qi)}
\,.
\tag 3.8
$$
This is equivalent to
$$
2 + 2 \left( \bold v_3 + \bold {\overline v}_3 \right) 
  \ \equiv \ 
  0 \pmod {2 (1 + \qi)}
\,.
\tag 3.9
$$
As the left side reduces to $4$, the equation is obviously valid.
\locqed
\enddemo

\bigskip
\head 4. The Correspondence Theorem \endhead

A quaternion algebra over a field is known to be isomorphic
to either a division algebra or the ring of $2 \times 2$
matrices over the field (see Pierce [7]).
Hurwitz's Correspondence Theorem is an analogue of this
fact for a certain quaternion order over the ring of rational
integers modulo an odd integer $m$.
In particular, the Correspondence theorem states there
is an isomorphism of the Hurwitz 
quaternions taken modulo odd $m$ with the $2 \times 2$
matrices over ${\Bbb Z} / m {\Bbb Z}$.
The Correspondence Theorem comes over in its entirety in
the case of the order $H_{1,2,2}$, as will be seen in this
section.

Let $m$ be an odd positive integer.
Since $m$ is in the centralizer, the left ideal, right ideal
and two sided ideal generated by $m$ are all the same.

\proclaim{Lemma 17}
Each quaternion $\bold q \in H_{1,2,2}$ is congruent modulo
$m$ to an element of $H_{1,2,2}^0$.
\endproclaim

\demo{Proof}
Set 
$\bold q = q_1 \bold v_1 + q_2 \bold v_2 + q_3 \bold v_3 + q_4 \bold v_4$.
Then 
$$
\bold q \ \equiv \ 
  q_1 \bold v_1 + q_2 \bold v_2 + (1 + m)\, q_3 \bold v_3 
  + (1 + m)\, q_4 \bold v_4
  \pmod m
\,.
\tag 4.1
$$
Since $1 + m$ is even $(1 + m) \bold v_3$ and $(1 + m) \bold v_4$
are in $H_{1,2,2}^0$.
\locqed
\enddemo

\proclaim{Lemma 18}
The $m^4$ quaternions
$$
\bold q \ = \ 
   q_1 \bold v_1 + q_2 \bold v_2 + q_3 \sqrt 2 \,\qj 
   + q_4 \sqrt 2 \,\qk \,,
     \quad\ 
   q_1,\, \dots,\, q_4 \,\in\, 
   \left\{  0,\, 1,\, \dots,\, m-1 \right\} 
\tag 4.2
$$
form a complete residue system for $H_{1,2,2} / (m)$.
\endproclaim

\demo{Proof}
By the previous Lemma, each element of $H_{1,2,2}$ is 
congruent modulo $m$ to a quaternion listed in (4.2).
Considered as free abelian groups, the quotient group, 
$H_{1,2,2} \,/\, (m)$ has order $m^4$ so we have the
proper number of elements.
\locqed
\enddemo

\proclaim{Lemma 19}
For an odd rational integer $m$, there exists 
$r,\, s \in \Bbb Z$ such that
$$
2^{-1} + r^2 + s^2 
    \ \equiv \ 
  0 \pmod m  
\tag 4.3
$$  
where $2^{-1}$ is the multiplicative inverse of 
$2$ modulo $m$.  
\endproclaim

\demo{Proof}
As in Hurwitz [6], it can be shown by elementary number
theory that for odd $m$ there exist $r,\, s \in \Bbb Z$
such that
$2 + r^2 + s^2 \equiv 0 \pmod m$.
Multiply this equation through by the square of $2^{-1}$
to prove the Lemma.
\locqed
\enddemo

At this point we wish to prove the analogue of Hurwitz's
Correspondence Theorem for the case of the quadratic form
$x^2 + y^2 + 2 z^2 + 2 w^2$.
It turns out to be advantageous to start with appropriate
versions of the Hurwitz variables 
$\qxi_1,\, \qxi_2,\, \qxi_3,\, \qxi_4$.
Using these we then produce the $\alpha,\, \beta,\, \gamma,\,
\delta$ and the coefficients
$q_1,\, q_2,\, q_3,\, q_4$ of the quaternion $\bold q$.
See Hurwitz [6, Vorlesung 8] for details in the case of
the form $x^2 + y^2 + z^2 + w^2$.

For the remainder of this section, $m$ will be an odd rational
integer.
We fix such an $m$, and choose $r$, $s$ as in Lemma 19.

\definition{Definition 20}
$$
\eqalign{
\qxi_1 \ &= \ 1 + r \sqrt 2 \, \qj + s \sqrt 2 \, \qk \,,
    \quad \hphantom{-i}
  \qxi_2 \ = \ \qi + s \sqrt 2 \, \qj - r \sqrt 2 \, \qk
     \cr  
\qxi_3 \ &= \ -i + s \sqrt 2 \, \qj - r \sqrt 2 \, \qk \,,
     \quad \hphantom{1}
   \qxi_4 \ = \ 1 - r \sqrt 2 \, \qj  - s \sqrt 2 \, \qk \,.
      \cr
}
\tag 4.4
$$
\enddefinition

\proclaim{Lemma 21}
The following orthogonal relations hold modulo $m$.
$$
\eqalign{
  &\qxi_1 \, \qxi_3 \ \equiv \ \qxi_1 \, \qxi_4
      \ \equiv \ 0 \,, \qquad
   \qxi_2 \, \qxi_1 \ \equiv \   \qxi_3 \,\qxi_4 
      \ \equiv \   0 \,,
     \cr    
  &\qxi_4 \, \qxi_1 \ \equiv \ \qxi_4 \, \qxi_2
       \ \equiv \ 0 \,, \qquad
   \qxi_2^2 \ \equiv \ \qxi_3^2 \ \equiv \ 0 \,.
      \cr     
}
\tag 4.5
$$
In addition, the non-orthogonal relations below also are valid modulo $m$.
$$
\eqalign{
 &\qxi_1^2 \ \equiv \ \qxi_2 \,\qxi_3 \ \equiv \  2 \,\qxi_1 \,, 
       \qquad \hphantom{\qxi_3}
    \qxi_1 \,\qxi_2 \ \equiv \ \qxi_2 \,\qxi_4  \ \equiv \  2\, \qxi_2 \,,
    \cr
 &\qxi_3 \,\qxi_1 \ \equiv \ \qxi_4 \,\qxi_3 \ \equiv \ 2 \,\qxi_3 \,, \qquad
     \qxi_3 \,\qxi_2 \ \equiv \ \qxi_4^2 \ \equiv \ 2 \,\qxi_4 \,.
\cr      
}
\tag 4.6
$$

\endproclaim

\demo{Proof}
Verified by computer algebra.
\locqed
\enddemo

\definition{Definition 22}
Given a four-tuple of integers 
$\alpha,\, \beta,\, \gamma,\, \delta$,
modulo odd $m$,
the corresponding quaternion $\bold q \in H_{1,2,2}/(m)$ 
is defined by
$$
\eqalign{
  2\, \bold q \ &= \ 
    2 \, q_1 + 2 \,q_2 \qi + 2 \, q_3 \sqrt 2 \, \qj
    + 2 \,q_4 \sqrt 2 \, \qk
  \cr
  \ &\equiv \   
    \alpha \,\qxi_1 + \beta \,\qxi_2 + \gamma \,\qxi_3
      + \delta \,\qxi_4
      \pmod m   \,.
}
\tag 4.7
$$
\enddefinition

It is clear from Definition 20 that the quantities 
$2\, q_1,\, \dots,\, 2\, q_4$ are well defined modulo $m$.
Since $m$ is odd that means $\bold q$ is well defined modulo
$m$.
Additionally we have the following.

\proclaim{Lemma 23}
With $\bold q$, and $\alpha$, $\beta$, $\gamma$,
and $\delta$ as above, we have modulo $m$
$$
\eqalign{
  2 \, q_1 \ &\equiv \ \alpha + \delta \,, \qquad
    2 \, q_3 \ \equiv \ r (\alpha - \delta)  + s (\beta + \gamma)
  \cr  
  2 \, q_2 \ &\equiv \ \beta - \gamma \,, \qquad
    2 \, q_4 \ \equiv \ s (\alpha - \delta) - r (\beta + \gamma)
  \cr  
}
\tag 4.8
$$
\endproclaim

\demo{Proof}
This is easily verified by hand computation or computer algebra.
\enddemo

We now produce the inverse mapping from quaternions 
$\bold q \in H_{1,2,2}/(m)$ to the four-tuple
that are coefficients of $\qxi_1$, $\qxi_2$, $\qxi_3$, $\qxi_4$.
This will later help demonstrate the one to one correspondence
between modulo $m$ quaternions and certain $2 \times 2$
integer matrices.

\proclaim{Lemma 24}
The inverse mapping from the quaternion $\bold q$ of\/ {\rm (4.7)}
to the modulo $m$ four-tuple
$\alpha$, $\beta$, $\gamma$, $\delta$ is given by
$$
\eqalign{
  \alpha \ &\equiv \ q_1 - 2 \,r \,q_3 - 2 \,s \,q_4 \,, \qquad
     \beta \ \equiv \ q_2 - 2 \,s \,q_3 + 2 \,r \,q_4 \,,
  \cr     
  \delta \ &\equiv \ q_1 + 2 \,r \,q_3 + 2 \,s \,q_4 \,, \qquad
     \gamma \ \equiv \ -q_2 - 2 \,s \,q_3 + 2 \,r \,q_4 \,,
  \cr     
}
\tag 4.9
$$
where all congruences are taken modulo $m$.
\endproclaim

\demo{Proof}
Given the relations of the previous Lemma we compute
$$
\eqalign{
  2 \,r \, q_3 + 2 \, s \, q_4
    \ &\equiv \ 
       r^2 (\alpha - \delta) + rs \,(\beta + \gamma) + s^2(\alpha - \delta)
          - sr \, (\beta + \gamma)
  \cr	        
  \ &\equiv \ 
     (r^2 + s^2) \, (\alpha - \delta)
     \ \equiv \ -2^{-1} (\alpha - \delta)
     \pmod m  \,.      
}
\tag 4.10
$$
Thus $-4\, r\, q_3 - 4\, s\, q_4 \,\equiv\, \alpha - \delta \pmod m$.
Using the relationship for $q_1$ we find
$$
\eqalign{
  &2\, q_1 - 4\, r\, q_3 - 4\, s\, q_4
    \ \equiv \ 
    2\, \alpha \pmod m
  \cr
  \Longrightarrow  \ \    
     &\alpha \ \equiv \
     q_1 - 2\, r\, q_3 - 2\, s\, q_4  \pmod m
  \cr
}
\tag 4.11
$$
Similarly we find
$\delta \equiv q_1 + 2\, r\, q_3 + 2\, s\, q_4 \pmod m$.

The formulas for $\beta$ and $\gamma$ are deduced in a similar
fashion.
\locqed
\enddemo

\proclaim{Lemma 25}
Let the quaternion $\bold q$ and $\alpha$, $\beta$, $\gamma$, $\delta$ 
be related as in the previous Lemma.
Then
$$
N ( \bold q) \ \equiv \ 
  \alpha \, \delta \,-\, \beta \, \gamma \pmod m
  \,.
\tag 4.12
$$
\endproclaim

\demo{Proof}
For the left side we have, modulo $m$,
$\bold q \equiv q_1 + q_2\, \qi + q_3 \sqrt 2\, \qj + q_4 \sqrt 2\,\qk$ 
which gives
$N( \bold q ) \equiv q_1^2 + q_2^2 + 2\, q_3^2 + 2\, q_4^2$.
On the other side of the equation
$$
\eqalign{
\alpha \, \delta \,-\, \beta \, \gamma 
    \ &\equiv \ 
  \left[  
          q_1^2 - (2\, r\, q_3 + 2\, s\, q_4)^2
  \right]
     \,-\,
  \left[
           (-2\, s\, q_3 + 2\, r\, q_4)^2 - q_2^2
  \right]
    \cr
    \ &\equiv \ 
  q_1^2 + q_2^2 - (4\, r^2 + 4\, s^2)\, q_3^2 
                - (4\, r^2 + 4\, s^2)\, q_4^2  
    \cr    
}
\tag 4.13
$$
but $r^2 + s^2 \equiv -2^{-1}$ so 
$\alpha \, \delta \,-\, \beta \, \gamma$
reduces to
$q_1^2 + q_2^2 + 2\, q_3^2 + 2\, q_4^2$ modulo $m$.
\locqed
\enddemo

\proclaim{Theorem 26 (Hurwitz Correspondence Theorem)}
There is a ring isomorphism from the quaternions of
$H_{1,2,2}$ modulo $m$ to 
the ring of $2 \times 2$ matrices of rational integers
modulo $m$,
$M_2 ({\Bbb Z} / m {\Bbb Z} )$.
Using the notation of the previous few Lemmas, the 
isomorphism is given using equation {\rm (4.9)}:
$$
\tau \,:\, \bold q
  \ \longrightarrow \ 
    \pmatrix
       \alpha & \beta  \cr
       \gamma & \delta \cr
    \endpmatrix
    \,.
\tag 4.14
$$
This isomorphism maps the quaternion norm into the determinant
of the corresponding matrix.
\endproclaim

\demo{Proof}
Since a complete residue system modulo $m$ for $H_{1,2,2}$ is
given by (4.2), both $H_{1,2,2} / (m)$ and 
the ring of $2 \times 2$ matrices over ${\Bbb Z} / m {\Bbb Z}$ 
contain $m^4$ elements.
Thus the equations (4.8) and (4.9) produce 
mappings between $H_{1,2,2} / (m)$ and 
$M_2 ({\Bbb Z} / m {\Bbb Z} )$.
Since the mapping corresponding to equation (4.9) returns to
the preimage of the map of (4.8), the former mapping is surjective.
As the sets in question have the same number of elements,
the mapping of (4.9) is a bijection.
The same holds for the mapping of (4.8).
Note $\tau (1 + 0\qi +0 \sqrt 2 \qj + 0 \sqrt 2 \qk)$ is the
identity matrix
$\left( 
   \smallmatrix 
   1 & 0 \cr
   0 & 1 \cr
   \endsmallmatrix
 \right)$.
The correspondence between the norm and determinant was demonstrated
in Lemma 25.
Suppose
$$
\bold q \ \equiv \ q_1 + q_2 \,\qi + q_3 \sqrt 2 \,\qj 
                   + q_4 \sqrt 2 \,\qk 
  \,, \ \quad
\bold q^\prime \ \equiv \ q_1^\prime + q_2^\prime \,\qi
    + q_3^\prime \sqrt 2 \,\qj + q_4^\prime \sqrt 2 \,\qk  
\tag 4.15
$$
and 
$$
\tau (\bold q )
  \ = \ 
    \pmatrix
       \alpha & \beta  \cr
       \gamma & \delta \cr
    \endpmatrix
    \,, \qquad
\tau (\bold q^\prime )
  \ = \ 
    \pmatrix
       \alpha^\prime & \beta^\prime  \cr
       \gamma^\prime & \delta^\prime \cr
    \endpmatrix
    \,.
\tag 4.16
$$
Then
$$
\tau (\bold q + \bold q^\prime )
  \ = \ 
  \tau \left(
          (q_1 + q_1^\prime) + (q_2 + q_2^\prime) \qi
	  + (q_3 + q_3^\prime) \sqrt 2 \,\qj
	  + (q_4 + q_4^\prime) \sqrt 2 \,\qk
       \right)
  .  
\tag 4.17
$$
The linearity of the equations (4.9) implies 
$\tau ( \bold q + \bold q^\prime) \ = \ 
     \tau (\bold q) + \tau (\bold q^\prime)$.
To find $\tau (\bold q \, \bold q^\prime)$ consider     
$$
\eqalign{
  2\, \bold q 
     \ &\equiv \   
       \alpha \,\qxi_1 + \beta \,\qxi_2 + \gamma \,\qxi_3
         + \delta \,\qxi_4
         \,,
  \cr
  2\, \bold q^\prime
     \ &\equiv \   
       \alpha^\prime \,\qxi_1 + \beta^\prime \,\qxi_2 
         + \gamma^\prime \,\qxi_3 + \delta^\prime \,\qxi_4
          \pmod m   \,.
}
\tag 4.18
$$
Multiplying and simplifying using the relations (4.5) and (4.6)
we find
$$
\eqalign{
4 \, \bold q \, \bold q^\prime
  \ &\equiv \ 
  2\,\left[ \,
       (\alpha\, \alpha^\prime + \beta\, \gamma^\prime) \,\qxi_1
       + (\alpha\, \beta^\prime + \beta\, \delta^\prime) \,\qxi_2
    \right.   
  \cr       
     &\qquad\,\left.
       + (\gamma\, \alpha^\prime + \delta\, \gamma^\prime) \,\qxi_3
       + (\gamma\, \beta^\prime + \delta\, \delta^\prime) \,\qxi_4 \,
    \right]   
    \pmod m \,.
  \cr    
}    
\tag 4.19
$$
Thus
$$
\tau (\bold q \, \bold q^\prime )
  \ = \ 
    \pmatrix
       \alpha & \beta  \cr
       \gamma & \delta \cr
    \endpmatrix
  \cdot
    \pmatrix
       \alpha^\prime & \beta^\prime  \cr
       \gamma^\prime & \delta^\prime \cr
    \endpmatrix
  \ = \ 
    \tau (\bold q) \, \tau (\bold q^\prime)
    \,,
\tag 4.20
$$
so $\tau$ is indeed a ring isomorphism.
\locqed
\enddemo

We slightly modify Hurwitz's definition of a quaternion primitive
to $m$ for the sake of clarity.

\definition{Definition 27}
A quaternion
$\bold g = g_1 \bold v_1 + g_2 \bold v_2 + g_3 \bold v_3 + g_4 \bold v_4 
  \in H_{1,2,2}$
is primitive to $m$ if
$\gcd (g_1,\, g_2,\, g_3,\, g_4,\, m) = 1$.
\enddefinition

Clearly any quaternion congruent modulo $m$ to $\bold g$ in $H_{1,2,2}$
is primitive to $m$ if and only if $\bold g$ is primitive to $m$.
Therefore primitivity to $m$ is well defined on $H_{1,2,2}/(m)$.

\proclaim{Lemma 28}
A quaternion
$\bold q = q_1 + q_2 \qi + q_3 \sqrt 2 \,\qj + q_4 \sqrt 2\,\qk \in
  H_{1,2,2}^0$
is primitive to $m$ if and only if
$\gcd (q_1,\, q_2,\, q_3,\, q_4,\, m) = 1$.
\endproclaim

\demo{Proof.}
In terms of the basis for $H_{1,2,2}$
$$
\bold q 
  \ = \  q_1 + q_2 \qi + q_3 \sqrt 2 \,\qj + q_4 \sqrt 2\,\qk 
  \ = \  g_1 \bold v_1 + g_2 \bold v_2 + g_3 \bold v_3 + g_4 \bold v_4 
  \,,
\tag 4.21
$$
thus it follows that
$$
\eqalign{
    &2\, q_1 \ = \ 2\, g_1 + g_3 + g_4\,, \qquad
    2\, q_3 \ = \ g_3\,,
  \cr
    &2\, q_2 \ = \ 2\, g_2 + g_3 + g_4\,, \qquad
    2\, q_4 \ = \ g_4 \,.
  \cr 
}
\tag 4.22
$$
For $m$ odd
$$
\eqalign{
  \gcd\, (g_1,\, g_2,\, g_3,\, g_4,\, m) 
      \ &= \ 
    \gcd\, (2\, g_1,\, 2\, g_2,\, g_3,\, g_4,\, m)
      \cr
      \ &= \ 
    \gcd\, (2\, g_1 + g_3 + g_4,\, 2\, g_2 + g_3 + g_4,\, g_3,\, g_4,\, m) 
      \cr
      \ &= \ 
    \gcd\, (2\, q_1,\, 2\, q_2,\, g_3,\, g_4,\, m) 
      \cr
      \ &= \ 
    \gcd\, (2\, q_1,\, 2\, q_2,\, 2\, q_3,\, 2\, q_4,\, m) 
      \cr
      \ &= \ 
    \gcd\, (q_1,\, q_2,\, q_3,\, q_4,\, m) \,.
      \cr
}
\tag 4.23
$$
The Lemma follows immediately.
The definition of primitivity to $m$ used here is
thus equivalent to Hurwitz's definition.
\locqed\enddemo

\definition{Definition 29}
An element
$\left( \smallmatrix
          \alpha & \beta  \cr
          \gamma & \delta \cr
       \endsmallmatrix \right) \in M_2 ({\Bbb Z} / m {\Bbb Z})$ 
is primitive to $m$ if the greatest common divisor
$\gcd (\alpha,\, \beta,\, \gamma,\, \delta,\, m) = 1$.
\enddefinition

\proclaim{Theorem 30}
For $m$ odd, a quaternion primitive to $m$ corresponds to a
matrix primitive to $m$, and conversely.
\endproclaim

\demo{Proof}
From the transformation formulas (4.8) and (4.9), this is clear.
\locqed
\enddemo

\proclaim{Theorem 31}
The number of incongruent quaternions in 
$H_{1,2,2} \,/\, (m)$ that are primitive to $m$ and satisfy
$N (\bold q) \equiv 0 \pmod m$ is
$$
\psi (m)  \ = \ 
m^3 \, \prod_{p \,\mid\, m} \, 
    \left( 1 - {1 \over{p^2}} \right) 
    \left( 1 + {1 \over p} \right)
\tag 4.24
$$
where $p$ runs through all primes dividing $m$.
\endproclaim

\demo{Proof}
By the Hurwitz Correspondence Theorem, we need to count the
same matrices as in Hurwitz [6, Vorlesung 8].
Thus we get the same result.
\locqed
\enddemo

\proclaim{Theorem 32}
The number of incongruent quaternions in 
$H_{1,2,2} \,/\, (m)$ that satisfy
$N (\bold q) \equiv 1 \pmod m$ is
$$
m^3 \, \prod_{p \,\mid\, m} \, 
    \left( 1 - {1 \over{p^2}} \right) 
\tag 4.25
$$
where $p$ runs through all primes dividing $m$.
\endproclaim

\demo{Proof}
Same proof as in the previous Theorem.
\locqed
\enddemo

\proclaim{Lemma 33}
Two distinct units of $H_{1,2,2}$ cannot be congruent to each
other modulo $m$ for odd $m > 1$.
\endproclaim

\demo{Proof}
Suppose not.
Then $\bold u_1 \equiv \bold u_2 \pmod m$ for units $\bold u_1$
and $\bold u_2$.
Let $\bold n = \bold u_1 \cdot \bold u_2^{-1}$.
Then $\bold n$ is also a unit of $H_{1,2,2}$ and 
$\bold n \equiv 1 \pmod m$.
Since $m \mid  (\bold n - 1)$ it follows that
$m^2 \mid   N( \bold n -1 )$.
Computation shows that the values of $N( \bold n - 1)$ are in
the set $\{ 0,\, 1,\, 2,\, 3,\, 4 \}$ as $\bold n$ ranges through
the units.
As none of these values are divisible by the square of an odd
number, the Lemma follows.
\locqed
\enddemo

\proclaim{Theorem 34}
Let $p$ be an odd prime and $\bold f$ any of the
$\psi (p) = (p^2 - 1) (p+1)$ quaternions of $H_{1,2,2}$
modulo $p$ which satisfy
$N( \bold f ) \equiv 0 \pmod p$
but have at least one component not divisible by $p$.
The number of distinct quaternion solutions $\bold x$ modulo $p$
which satisfy $\bold x \, \bold f \equiv 0 \pmod p$ is $p^2$.
\endproclaim

\demo{Proof}
The proof is the same as in Theorem 31.
\locqed
\enddemo

\bigskip
\head 5. Prime quaternions \endhead

We note that $H_{1,2,2}$ is a division algebra as it is a
noncommutative subring of the set of all quaternions.
It is easy to see that for any nonzero quaternion, the
left multiplicative inverse and the right multiplicative
inverse are equal.
We recall that an element of a ring is a unit if it has
an inverse in the ring.
The units of $H_{1,2,2}$ are the elements of norm $1$ and
are listed in (1.3).

\definition{Definition 35}
A prime of $H_{1,2,2}$ is a nonzero nonunit $\qpi$ such that
if
$\qpi = \bold a \cdot \bold b$
in $H_{1,2,2}$ then at least one of $\bold a$ or $\bold b$
is a unit.
\enddefinition

\proclaim{Theorem 36}
A rational prime $p \in \Bbb Z$ is not a prime of $H_{1,2,2}$.
\endproclaim

\demo{Proof}
Clearly $2$ is not prime as $2 = (1 + \qi) \, (1 - \qi)$.
Consider the case of $p$ an odd rational prime.
By previous results we may choose $\bold q \in H_{1,2,2}$
with $\bold q$ primitive to $p$ and $N( \bold q) \equiv 0 \pmod p$.
Thus $\bold q \not\equiv 0 \pmod p$.

Set $\bold d$ equal to the right greatest common divisor of
$p$ and $\bold q$.
Thus there exists $\bold b,\, \bold c,\, \bold d_1,\, \bold d_2$
in $H_{1,2,2}$ such that 
$$
\bold d \ = \ 
    \bold c \, p + \bold b \, \bold q \,,\quad
  p \ = \ \bold d_1 \, \bold d \,, \quad
  \bold q \ = \ \bold d_2 \, \bold d 
\,.
\tag 5.1
$$
Taking conjugates we have 
$\bold {\overline d} = p \, \bold {\overline c} 
           + \bold {\overline q} \,\bold {\overline b}$ so
$$
\eqalign{
  \bold d \, \bold {\overline d}
    \ &= \ 
      p^2 \, \bold c \, \bold {\overline c}
      + p \,\bold c \,\bold {\overline q} \,\bold {\overline b}
      + p \,\bold b \,\bold q \,\bold {\overline c}
      + \bold b \,\bold q \,\bold {\overline q} \,\bold {\overline b}
  \cr
    \ &= \   
      p \,\bold k + N(\bold q) \, N( \bold b) 
  \cr
}
\tag 5.2
$$
for some $\bold k \in H_{1,2,2}$.
Thus $N(\bold d) = \bold d \,\bold {\overline d} \equiv 0 \pmod p$.
If $\bold d$ were a unit, then $N(\bold d) = 1 \not\equiv 0 \pmod p$.
Hence $\bold d$ is not a unit.

If $\bold d_1$ were a unit, then from $p = \bold d_1 \, \bold d$ it
follows that $\bold d = \bold d_1^{-1} \, p$.
Hence 
$\bold q = \bold d_2 \, \bold d = \bold d_2 \, \bold d_1^{-1} \, p$.
Thus $\bold q \equiv 0 \pmod p$ in $H_{1,2,2}$ so $\bold q$ is not
primitive to $p$, a contradiction.
Thus $p = \bold d_1 \, \bold d$ is a product of two elements neither
of which is a unit, so $p$ cannot be prime in $H_{1,2,2}$.
\locqed
\enddemo

\proclaim{Theorem 37}
An element $\qpi \in H_{1,2,2}$ is prime iff $N(\qpi)$
is a rational prime.
\endproclaim

\demo{Proof}
Suppose $N( \qpi ) = p$, a prime of $\Bbb Z$, and that $\qpi$
factors as $\bold a \, \bold b$ in $H_{1,2,2}$.
Then $p = N(\bold a ) \, N(\bold b)$ in $\Bbb Z$ so one of
$N(\bold a)$ or $N(\bold b)$ is $1$.
Hence one of $\bold a$ or $\bold b$ is a unit, so $\qpi$ is
a prime.

In the other direction, suppose $\qpi \in H_{1,2,2}$ is prime.
Since $\qpi$ is not a unit, $N(\qpi )$ is not $1$ so there exists
some rational prime $p$ dividing $N(\qpi)$.
Set $\bold d$ to the one sided greatest common divisor of $p$
and $\qpi$ as in the previous theorem.
As above, we find that $N(\bold d) \equiv 0 \pmod p$, hence
$\bold d$ cannot be a unit.
Also we have
$$
  p \ = \ \bold d_1 \, \bold d \,, \quad
  \qpi \ = \ \bold d_2 \, \bold d 
\tag 5.3
$$
for $\bold d_1,\, \bold d_2 \in H_{1,2,2}$.
Since $\qpi$ is prime, $\bold d_2$ must be a unit.
Thus $p = \bold d_1 \, \bold d_2^{-1} \, \qpi$.
Taking the norms, 
$p^2 = N(\bold d_1) \cdot 1 \cdot N(\qpi )$.

Since $p$ divides $N (\qpi)$, that means $N(\bold d_1)$ 
can be either $1$ or $p$.
If it were $1$ then $\bold d_1$ would be a unit implying that
$p$ is a unit times the prime $\qpi$.
Thus $p$ would be a prime of $H_{1,2,2}$, a contradiction.
We conclude that $N(\bold d_1) = p$ and consequently that
$N(\qpi ) = p$.
\locqed
\enddemo

We now wish to find all prime quaternions dividing a
prime rational number.

\proclaim{Lemma 38}
Let $\qpi$ be a prime rational quaternion dividing
$2$ in $H_{1,2,2}$.
Then $\qpi$ equals $1 + \qi$ times a unit.
\endproclaim

\demo{Proof}
The demonstration is the same as in Hurwitz [6, Vorlesung 9].
\locqed
\enddemo

\proclaim{Lemma 39}
Let $\bold f \in H_{1,2,2}$ be a quaternion primitive to a rational
prime $p \not= 2$.
Suppose $N(\bold f ) \equiv 0 \pmod p$.
Then there exists $\bold {\widehat f} \in H_{1,2,2}$ such that
$\bold {\widehat f} \equiv \bold f \pmod p$ and
$N( \bold {\widehat f} ) \equiv 0 \pmod p$ while
$N( \bold {\widehat f} ) \not\equiv 0 \pmod {p^2$}.
\endproclaim

\demo{Proof}
The proof is similar to that in Hurwitz [6, Vorlesung 9]
though the norm form is more complicated.
In this case we have
$$
\eqalign{
N ( a_1 \bold v_1 + a_2 \bold v_2 + a_3 \bold v_3
    + a_4 \bold v_4 )
  \ = \ 
    &a_1^2 + a_2^2 + a_3^2 + a_4^2 + a_1 a_3 + a_2 a_3
  \cr
      &+ a_1 a_4 + a_2 a_4 + a_3 a_4
    \,.
  \cr
}  
\tag 5.4
$$      
Set
$ \bold f = f_1 \bold v_1 + f_2 \bold v_2 + f_3 \bold v_3 + f_4 \bold v_4$
with $f_1,\,f_2,\, f_3,\, f_4 \in \Bbb Z$
and 
$\bold {\widehat f} = f + 
     p ( t_1 \bold v_1 + t_2 \bold v_2 + t_3 \bold v_3 + t_4 \bold v_4 )$
with $t_1,\, t_2,\, t_3,\, t_4 \in \Bbb Z$.
Computer algebra shows
$$
\eqalign{
N( \bold{\widehat f} ) \ = \ 
  N ( \bold f ) 
    + p \,( 
       &2\, f_4 t_4 + f_3 t_4 + f_2 t_4 + f_1 t_4 + f_4 t_3
  \cr
       &+ 2\, f_3 t_3  + f_2 t_3 + f_1 t_3 + f_4 t_2 + f_3 t_2
  \cr
      &+ 2\, f_2 t_2 + f_4 t_1 + f_3 t_1 + 2\, f_1 t_1
    )   
    \pmod {p^2}  \,.
  \cr  
}
\tag 5.5
$$     
We need only solve for rational integer 
$t_1,\, t_2,\, t_3,\, t_4$ such that
$$
\eqalign{
  &(f_4 + f_3 + 2\, f_1) \, t_1 + (f_4 + f_3 + 2\, f_2 )\, t_2
     \cr
  &+ (f_4 + 2\, f_3 + f_2 + f_1 )\, t_3
  + (2\, f_4 + f_3 + f_2 + f_1)\, t_4
  \ \equiv \ 1 \pmod p
  \,.
    \cr
}
\tag 5.6
$$     
If any of the coefficients of $t_1,\, t_2,\, t_3,\, t_4$ are not
congruent to zero modulo $p$ then set the other $t$'s to zero and
solve for the value of the remaining $t$ modulo $p$.
This produces an $\bold {\widehat f}$ satisfying the condition
of the Lemma.

The remaining possibility is if all four of the coefficients of
the $t$'s are congruent to zero modulo $p$.
From the coefficients of $t_1$ and $t_2$ we find
$$
2\, f_1 \ \equiv \ 2\, f_2 \ \equiv \ 
  -(f_3 + f_4) \pmod p
  \,.
\tag 5.7
$$
Adding the coefficients of $t_3$ and $t_4$ 
$$
3 ( f_3 + f_4) + 2 f_2 + 2 f_1
  \ \equiv \ 0 \pmod p
  \,.
\tag 5.8
$$     
Combining (5.7) with (5.8) results in
$f_3 + f_4 \equiv 0 \pmod p$.
This implies $f_1 \equiv f_2 \equiv 0 \pmod p$.
The coefficient of $t_3$ now reduces to
$f_4 + 2\, f_3$.
This being congruent to $0$ by hypothesis, implies
that $f_3$, and then $f_4$ are congruent to $0$ modulo
$p$.
Hence $f$ is not primitive to $p$, a contradiction.
\locqed
\enddemo

\proclaim{Lemma 40}
Let $\bold f \in H_{1,2,2}$ be a quaternion primitive to a rational
prime $p \not= 2$.
Suppose $N(\bold f ) \equiv 0 \pmod p$.
Then we can uniquely associate $\bold f$ to a primary prime quaternion
$\qpi$ of norm $p$ by the following algorithm.

Find\/ $\bold {\widehat f} \in H_{1,2,2}$ such that\/
$\bold {\widehat f} \equiv \bold f \pmod p$ and
$N( \bold {\widehat f} ) \equiv 0 \pmod p$ while
$N( \bold {\widehat f} ) \not\equiv 0 \pmod {p^2$}.
Compute the right greatest common divisor of\/ $\bold{\widehat f}$
and $p$.
This is a prime quaternion of norm $p$.
Let $\qpi$ be the primary associate of this prime quaternion.
\endproclaim

\demo{Proof}
We must show $\qpi$ exists with the requisite properties and
is unique.
$\bold{\widehat f}$ exists by the Lemma above.
Let $\qtau$ be the right greatest common divisor of 
$\bold{\widehat f}$ and $p$.
By an argument very similar to (5.1) and (5.2) we find
$N(\qtau ) \equiv 0 \pmod p$.
Hence $\qtau$ cannot be a unit.
Additionally $N(\qtau)$ divides $N(p) = p^2$, so $\qtau$ is odd.
We may take $\qpi$ to be the unique primary associate of $\qtau$.

We have $\bold {\widehat f} = \bold a \, \qtau$ for some
$\bold a \in H_{1,2,2}$.
However
$N(\bold a) \, N(\qtau) \equiv N( \bold{\widehat f}) 
   \not\equiv 0 \pmod {p^2}$.
Hence $N(\qtau)$ cannot equal $0$ modulo $p^2$.
Thus $\qtau$ is of norm $p$, and must be a prime.
Since $N(\qpi) = N(\qtau)$ it follows that $\qpi$ is also
prime, and indeed a primary prime of norm $p$.

We must show that the $\qpi$ constructed above is unique.
Let $\bold g \equiv \bold {\widehat g} \pmod p$ with
$N(\bold g) \equiv N(\bold {\widehat g}) \equiv 0 \pmod p$
while both of $N(\bold g)$ and $N(\bold {\widehat g})$ are
not congruent to $0$ modulo $p^2$.
Let the right greatest common divisor mapping described above
produce $\qpi$ for $\bold g$ and $\qtau$ for $\bold{\widehat g}$
where $\qpi$ and $\qtau$ are both primary primes of norm $p$.
Then $\bold g  =  \bold a \, \qpi$,
$\bold {\widehat g}  =  \bold b \, \qtau$ for some
$\bold a$, $\bold b$ in $H_{1,2,2}$.

From $N(\bold g)= N(\bold a) N(\qpi)$ it follows that
$N(\bold a) \not\equiv 0 \pmod p$.
If either one sided greatest common divisor of $\bold a$ and 
$p$ had norm greater than $1$, then $\bold a$ would have a
divisor of norm $p$ or $p^2$.
This would imply that $N(\bold a) \equiv 0 \pmod p$, a
contradiction.
Hence each one sided greatest common divisor of $\bold a$
and $p$ is a unit.
In particular there exists $\bold w,\, \bold x \in H_{1,2,2}$
such that $\bold w \,\bold a + \bold x \, p = 1$.
Thus $\bold w \,\bold a \,\qpi + \bold x \,p \,\qpi = \qpi$.
This implies $\bold w \,\bold g + \bold x \,\qpi \,p = \qpi$.
By hypothesis, $\bold g = \bold{\widehat g} + \bold c \, p$ for
some $\bold c \in H_{1,2,2}$.
Thus
$$
\qpi \ = \ 
  \bold w \,\bold{\widehat g} + \bold w \,\bold c \,p
    + \bold x \,\qpi p
  \ = \ 
  \bold w \,\bold b \,\qtau + (\bold w \,\bold c \,\overline\qtau
     + \bold x \,\qpi \,\overline\qtau ) \,\qtau  
  \,.
\tag 5.9
$$    
Thus $\qtau$ divides $\qpi$, but as both are primary primes,
they are equal.
\locqed
\enddemo

\proclaim{Lemma 41}
The mapping of the previous Lemma is surjective.
\endproclaim

\demo{Proof}
Let $\qpi$ be a primary prime quaternion of norm $p$.
Then we may take $\bold f = \qpi$ as $\qpi$ is primitive
to $p$ (else the norm of $\qpi$ would be divisible by $p^2$).
Since $N(\qpi) = p \not\equiv 0 \pmod {p^2}$ we may
take $\bold {\widehat f} = \qpi$.
As $p = N(\qpi) = \overline\qpi \, \qpi$, when we take the
primary form of the right greatest common divisor of $p$ and 
$\qpi$ we merely reproduce $\qpi$.
\locqed
\enddemo

\proclaim{Lemma 42}
A prime $\qpi$ divides a rational prime $p$ iff $N(\qpi) = p$.
\endproclaim

\demo{Proof}
$N(\qpi) = p$ implies 
$p = \overline\qpi \, \qpi = \qpi \, \overline\qpi$.
Thus $\qpi$ divides $p$ on the left and the right.

In the other direction, if $\qpi$ divides a rational prime $p$
then $N(\qpi) \mid N(p)$ so $N(\qpi)$ divides $p^2$.
Since $\qpi$ is prime, its norm must be a rational prime, and
hence equals $p$.
\locqed
\enddemo

\proclaim{Lemma 43}
Under the conditions of Lemma 40
two quaternions primitive to $p$, $\bold f$ and $\qfone$,
produce the same primary prime quaternion under the right 
greatest common divisor algorithm iff there exists
$\bold q \in H_{1,2,2}$ such that
$\qfone \equiv \bold q \, \bold f \pmod p$.
\endproclaim

\demo{Proof}
By previous results we may construct quaternions in $H_{1,2,2}$
congruent to each of $\bold f$ and $\qfone$ modulo $p$
such that the norms are divisible by $p$ but not by $p^2$.
We may replace $\bold f$ and $\qfone$ by the constructed quaternions.
and it is clear that the Lemma will be valid for the original 
quaternions iff it is valid for the constructed quaternions.

Suppose $\bold f$ and $\qfone$ produce the same quaternion $\qpi$
under the right gcd algorithm.
Then there exists $\qalpha,\, \qalfone \in H_{1,2,2}$ such that
$$
\bold f \ = \ \qalpha \, \qpi \,,\qquad
  \qfone \ = \ \qalfone \, \qpi
\tag 5.10
$$
with $N(\qpi) = p$.
Thus $N(\qalpha)$ and $N(\qalfone)$ are relatively prime to $p$
since otherwise $N(\bold f)$ or $N(\qfone)$ would be divisible by
$p^2$.
Thus we may solve the following congruence for the unknown
quaternion $\bold q$ 
$$
\bold q \, N( \qalpha ) \ \equiv \ 
  \qalfone \, \overline \qalpha  \pmod p 
\tag 5.11
$$
by merely considering coefficients of basis elements of $H_{1,2,2}$.
Since $N(\qalpha) = \qalpha \, \overline\qalpha$ we find
$$
\left(
  \bold q \, \qalpha - \qalfone \right) \, \overline\qalpha
    \ \equiv \ 
    0 \pmod p \,.
\tag 5.12
$$
Consider either one sided gcd of $\overline\qalpha$ and $p$.
If this gcd were not a unit then $p$ would divide $N(\overline\qalpha)$
and thus $N(\qalpha)$, a contradiction.
Thus each one sided gcd is a unit and thus there exists one sided linear
combinations of $\overline\qalpha$ and $p$ that sum up to $1$.
In particular there exists $\bold b \in H_{1,2,2}$ such that
$\overline\qalpha \, \bold b \equiv 1 \pmod p$.
Multiplying (5.12) by $\bold b$ on the right results in
$\bold q \, \qalpha - \qalfone \equiv  0 \pmod p$.
Multiplying by $\qpi$ on the right results in
$\bold q \, \qalpha \,\qpi - \qalfone \,\qpi \equiv  0 \pmod p$.
Thus 
$\bold q \, \bold f \equiv \qfone  \pmod p$.

In the other direction suppose that there exists $\bold q \in H_{1,2,2}$
such that
$\qfone \equiv \bold q \, \bold f \pmod p$
and that $\bold f$ corresponds to $\qpi$ under the right gcd algorithm.
Then there exists $\qalpha \in H_{1,2,2}$ such that
$\bold f = \qalpha \, \qpi$ and $N(\qpi) = p$.
For some $\bold c \in H_{1,2,2}$
$$
\qfone \ = \ \bold q \, \bold f + \bold c \, p
  \ = \ 
     \bold q \, \qalpha \,\qpi  + \bold c \, \overline\qpi \,\qpi
  \,.
\tag 5.13
$$
Thus $\qpi$ divides the right greatest common divisor of $\qfone$
and $p$.
Yet this right gcd has norm $p$ and the algorithm chooses the
primary associate.
The primary associate must therefore be $\qpi$ and hence
$\qfone$ also corresponds to $\qpi$.     
\locqed\enddemo

\proclaim{Lemma 44}
Let $\qpi$ be a primary prime quaternion of norm $p \not= 2$ 
in $H_{1,2,2}$.
The number of incongruent quaternions 
in $H_{1,2,2}/(p)$
primitive to $p$ that produce $\qpi$
under the right gcd algorithm is $p^2-1$.
\endproclaim

\demo{Proof}
The demonstration is the same as in Hurwitz [6, Vorlesung 9].
\locqed\enddemo

\proclaim{Corollary 45}
The number of primary prime quaternions in $H_{1,2,2}$ of norm
$p \not = 2$ is $p+1$.
\endproclaim

\demo{Proof}
The demonstration is the same as in Hurwitz [6, Vorlesung 9].
\locqed\enddemo

\proclaim{Corollary 46}
There are exactly $24 \, (p+1)$ prime quaternions dividing
an odd rational prime $p \not= 2$ in $H_{1,2,2}$.
\endproclaim

\demo{Proof}
The demonstration is the same as in Hurwitz [6, Vorlesung 9].
\locqed\enddemo

\bigskip
\head 6. Decomposition Theorem in $H_{1,2,2}$   \endhead

We wish to study the factorization into prime quaternions
in $H_{1,2,2}$.
We also would like to discover the number of primary
quaternions of given norm $m$ in this quaternion order.
As in Hurwitz [6, Vorlesung 10] any $\bold a \in H_{1,2,2}$
may be written as $\bold a = (1 + \qi)^r \, \bold b$ with
nonnegative integer $r$ and $\bold b$ an odd quaternion.
We note that $\bold b$ is only a unit factor from an odd
primary quaternion.

\definition{Definition 47}
A quaternion 
$\bold c = c_1 \bold v_1 + c_2 \bold v_2 + c_3 \bold v_3 + c_4 \bold v_4$
in $H_{1,2,2}$ is primitive if it is primary and
$\gcd (c_1,\, c_2,\, c_3,\, c_4) = 1$ in $\Bbb Z$.
\enddefinition

\proclaim{Lemma 48}
Let $\bold b$ be a primary quaternion and let $m$ equal the
greatest common divisor of the components of $\bold b$ with 
respect to the basis
$\{ \bold v_1,\, \bold v_2,\, \bold v_3,\, \bold v_4 \}$.
There exists a primitive quaternion $\bold c$ such that
$\bold b = m \,\bold c$ or $\bold b = -m \, \bold c$.
\endproclaim

\demo{Proof}
Since $\bold b$ is primary, $\bold b \equiv 1$ or $1 + 2\bold v_3$
modulo $2 (1 + \qi)$.
Thus $\bold b \equiv 1 \pmod 2$ which implies 
$\bold b \equiv 1 \pmod {1 + \qi}$.
Thus $\bold b$ is odd and hence $m$ is odd.
The proof now proceeds as in Hurwitz [6, Vorlesung 10].
\locqed\enddemo

It is simple to observe that a primary prime quaternion is also 
primitive.
While Hurwitz redefines conjugate at this point, we introduce
the definition of a variation of conjugate.

\definition{Definition 49}
Let $\qpi$ be a primary prime quaternion.
The {\sl p}-conjugate of $\qpi$, denoted $\overline\qpi_p$, is
$$
\overline\qpi_p \ = \ 
    \cases
           \overline\qpi & \hbox{if } \qpi \equiv 1 \pmod {2(1 + \qi)} \cr
           -\overline\qpi & \hbox{if } \qpi \equiv 1 + 2 \bold v_3
	                    \pmod {2(1 + \qi)} \,. \cr          
    \endcases
\tag 6.1
$$
\enddefinition

\proclaim{Lemma 50}
Let $\qpi$ be a primary prime quaternion.
The {\sl p}-conjugate of $\qpi$ is $\overline\qpi$ or
$-\overline\qpi$ depending on whether 
$p = N( \qpi ) \equiv 1$ or $-1 \pmod 4$ respectively.
\endproclaim

\demo{Proof}
We first note that if $k \in \Bbb Z$ is divisible by 
$1 + \qi$ then $k$ is even, for if $k$ were odd then
$N(k)$ would be odd but divisible by $N(1 + \qi) = 2$,
a contradiction.

If $\qpi \equiv 1 \pmod {2(1 + \qi)}$ then
$\overline\qpi \equiv 1 \pmod {2(1 + \qi)}$ so
$p = \qpi \,\overline\qpi \equiv 1 \pmod {2(1 + \qi)}$.
Thus $(p-1)/2$ is a rational integer divisible by $1 + \qi$
and hence even, so $p \equiv 1 \pmod 4$.

If $\qpi \equiv 1 + 2 \bold v_3 \pmod {2(1+\qi)}$ then using
equation (3.8)
$$
-\overline\qpi \ \equiv \ 
  - (1 + 2 \bold {\overline v}_3)  \ \equiv \ 
  1 + 2 \bold v_3 \pmod  {2(1 + \qi)}
\,.
\tag 6.2
$$
Thus 
$-p = - \qpi \, \overline\qpi \equiv (1 + 2 \bold v_3)^2
   \equiv 1 \pmod {2(1 + \qi)}$
so $(-p-1) / 2$ is a rational integer divisible by $1 + \qi$.
We conclude it is even, so $p \equiv -1 \pmod 4$.
\locqed\enddemo

Note that the {\sl p}-conjugate of a primary prime quaternion
is also primary by equation (3.8).

\proclaim{Lemma 51}
Let $\bold c$ be a primitive quaternion and suppose the rational
prime $p$ divides $N(\bold c )$.
Then the left greatest common divisor of $\bold c$ and $p$
is a prime quaternion $\qpi \in H_{1,2,2}$.
If chosen as primary, $\qpi$ is uniquely determined.
\endproclaim

\demo{Proof}
Let $\qlambda$ be the left greatest common divisor of $\bold c$
and $p$.
Then $N(\qlambda)$ divides $N(p)$ in $\Bbb Z$ so $N(\qlambda)$
can take only the values $1,\, p$ or $p^2$.

Suppose $N(\qlambda) = p^2$.
Let $p = \qlambda \, \qtau$ for some $\qtau \in H_{1,2,2}$.
Then $N(\qtau) = 1$ so $\qtau$ is a unit.
Thus $\qlambda$ is $p$ times a unit, which implies $p$ divides
$\bold c$ in $H_{1,2,2}$.
So $\bold c = p (g_1 \bold v_1 + g_2 \bold v_2 + g_3 \bold v_3
     + g_4 \bold v_4)$
in $H_{1,2,2}$.
We observe that $p$ divides the greatest common divisor of the
components of $\bold c$ which implies $\bold c$ is not primitive,
a contradiction.

Suppose $N(\qlambda ) = 1$.
Then $\qlambda$ is a unit.     
There exists $\bold a,\, \bold b \in H_{1,2,2}$ such that
$1 = \bold c \, \bold a + p \, \bold b$ from the left greatest
common divisor hypothesis.
Thus $1 \equiv \bold c \, \bold a \pmod p$.
Taking conjugates we find 
$1 \equiv \overline{\bold a} \, \overline{\bold c} \pmod p$.
Multiplying together $1 \equiv N(\bold c) N(\bold a) \pmod p$.
Thus $p$ does not divide $N(\bold c)$, a contradiction.

The only possibility remaining is $N(\qlambda) = p$ so 
$\qlambda$ is a prime quaternion.
\locqed\enddemo

\proclaim{Theorem 52 (Decomposition Theorem)}
Let $\bold c$ be a primitive quaternion of $H_{1,2,2}$ whose 
norm has the rational prime decomposition
$N(\bold c ) = p_1 \,p_2 \,p_3 \,\cdots$.
Then $\bold c$ can be uniquely represented in the form
$\bold c = \qpi_1 \,\qpi_2 \,\qpi_3 \,\cdots$
where $\qpi_1,\, \qpi_2,\, \dots$ are primary prime quaternions
of norm $p_1,\, p_2,\, \dots$ respectively.
Furthermore, this can be done for any rearrangement of the order
of the primes $p_1,\, p_2,\, \dots\,\,$.
\endproclaim

\demo{Proof}
The proof follows as in Hurwitz [6, Vorlesung 10].
\locqed\enddemo

\proclaim{Corollary 53}
Let $\bold c$ be a primitive quaternion of $H_{1,2,2}$ 
whose norm has the rational prime decomposition
$N(\bold c ) = p^h \,q^k \,\cdots$.
where $p,\, q\,\ \dots$ are distinct primes.
Then $\bold c$ can be uniquely factored into the form
$$
\bold c \ = \ 
  \qpi_1 \, \qpi_2 \,\cdots\, \qpi_h  \,\,
  \qchi_1 \, \qchi_1 \,\cdots\, \qchi_k \,\,
  \cdots
\tag 6.3
$$
where the $\qpi$'s are primary prime quaternions of norm $p$,
the $\qchi$'s are primary prime quaternions of norm $q$, etc.
\endproclaim

\demo{Proof}
The proof follows, as in Hurwitz [6, Vorlesung 10],
from the Decomposition Theorem.
\locqed\enddemo

We note that in the primary prime decomposition of $\bold c$
as above, a prime $\qpi$ and its {\sl p}-conjugate cannot appear
adjacent to each other.
If they were to do so, then $p = N(\qpi)$ would divide $\bold c$,
and thus all the components of $\bold c$ written in the basis
$\{ \bold v_1,\, \bold v_2,\, \bold v_3,\, \bold v_4 \}$.
Thus $\bold c$ would not be primitive.

\proclaim{Theorem 54}
Let $p,\, q\,\ \dots$ be rational primes and
$\qpi_s,\, \qchi_s,\, \dots$ be primary prime quaternions
in $H_{1,2,2}$ of norm $p,\, q,\, \dots$ respectively for
all positive integers $s$.
Then the product
$\qpi_1 \,\qpi_2 \,\cdots \,\qpi_h \,\,
  \qchi_1 \,\qchi_2 \,\cdots \,\qchi_k \,\cdots$
is a primitive quaternion as long as two {\sl p}-conjugate prime
quaternions never appear adjacent to each other.
\endproclaim

\demo{Proof}
The proof is the same as in Hurwitz [6, Vorlesung 10].
\locqed\enddemo

\proclaim{Theorem 55}
Let $m$ be an odd rational integer with the decomposition into
distinct primes
$m = p_1^{h_1} \, p_2^{h_2} \, \cdots$.
Then the number of primitive quaternions of norm $m$,
denoted $Q(m)$, is
$$
Q (m) \ = \ 
   m \,\prod_s \, \left(  1 + {1 \over {p_s}} \right) 
\tag 6.4
$$
\endproclaim

\demo{Proof}
The proof is the same as in Hurwitz [6, Vorlesung 10].
\locqed\enddemo

For computational reasons we define $Q(1) = 1$.
It is easy to see that for relatively prime odd rational integers
$m_1,\, m_2$ we have
$Q (m_1\, m_2 ) = Q(m_1) \, Q(m_2)$.

\proclaim{Theorem 56}
Let $m$ be an odd rational integer.
The number of primary quaternions in $H_{1,2,2}$ of norm $m$
is $\sigma (m)$, the sum of the divisors of $m$.
\endproclaim

\demo{Proof}
The proof is the same as in Hurwitz [6, Vorlesung 10].
\locqed\enddemo

\bigskip
\head 7. Counting solutions     \endhead

We have now developed enough number theory in the order
$H_{1,2,2}$ to determine the number of representations
in $\Bbb Z$ by the norm form.
As the expansions of the units of this order in the 
standard quaternion basis are relevant
to counting representations, it is convenient to list them,
namely (see Deutsch [1])
$$
\pm \left\{
  1,\,\,\, \qi,\,\,\,
  {1 \over 2} \pm {1 \over 2} \qi \pm {{\sqrt 2} \over 2} \qj,\,\,\,
  {1 \over 2} \pm {1 \over 2} \qi \pm {{\sqrt 2} \over 2} \qk,\,\,\,
  {{\sqrt 2} \over 2} \qj \pm {{\sqrt 2} \over 2} \qk
\right\}
\,.
\tag 7.1
$$

\proclaim{Theorem 57}
Let $n$ be a positive integer, $n = 2^r \, m$ with $m$ an odd
rational integer.
The number of representations of $n$ by the quadratic form
$x^2 + y^2 + 2\, z^2 + 2\, w^2$ is
$$
\cases
   4\, \sigma (m)  & \hbox{ if $n$ is odd} \cr
   8\, \sigma (m)  & \hbox{ if } r \,=\, 1 \cr
  24\, \sigma (m)  & \hbox{ if } r \,\geq\, 2 \cr
\endcases
\tag 7.2
$$
\endproclaim

\demo{Proof}
We wish to find the number of quaternions $\bold x \in H_{1,2,2}^0$
such that $N(\bold x) = n$.
Since $n = 2^r\, m$ we may write $\bold x = (1 + \qi)^r \, \bold y$
where $\bold y$ is an odd quaternion in $H_{1,2,2}$.
Being odd, $\bold y$ can be uniquely be factored as a unit times
a primary quaternion.
Let $\bold y = \bold u \, \bold c$ where $\bold u$ is a unit and
$\bold c$ is primary.
Then $\bold c \equiv 1 \pmod 2$ so there exists $\bold g \in H_{1,2,2}$
such that $\bold c = 1 + 2\, \bold g$.
Thus
$$
\eqalign{
  \bold x \ &= \ (1 + \qi)^r \, \bold u \, \bold c
             = \ (1 + \qi)^r \, \bold u \, (1 + 2\, \bold g)
      \cr 
  &= \ (1 + \qi)^r \, \bold u 
   + 2\, (1 + \qi)^r \, \bold u \, \bold g
   \,.
      \cr 
}
\tag 7.3
$$
Since $(1 + \qi)^r \, \bold u \, \bold g \in H_{1,2,2}$,
twice this quantity is in $H_{1,2,2}^0$.
Thus $\bold x \in H_{1,2,2}^0$ if and only if
$(1 + \qi)^r \, \bold u \in H_{1,2,2}^0$.
This is amenable to computation as there are only 24 units.

For $r = 0$ the only units in $H_{1,2,2}^0$ are
$\pm \{ 1,\, \qi \}$.
Hence the number of solutions $\bold x$ is four times the number
of primary quaternions of norm $m$, i.e.~$4 \, \sigma (m)$.

For $r = 1$ computation shows that the only units $\bold u$ such
that $(1 + \qi) \, \bold u$ are in $H_{1,2,2}^0$ are
$$
\pm \left\{
  1,\  \qi,\  
  {{\sqrt 2} \over 2} \qj \pm {{\sqrt 2} \over 2} \qk
  \right\}
\,.
\tag 7.4
$$
Thus the number of solutions $\bold x$ is $8 \, \sigma (m)$.

For $r \geq 2$ note that $(1 + \qi)^r$ contains a factor of $2$.
Hence $(1 + \qi)^r \, \bold u$ is in $H_{1,2,2}^0$ for all units
$\bold u$.
The number of solutions $\bold x$ is therefore 
$24 \, \sigma (m)$.
\locqed\enddemo

In the above Theorem we have cases where only $4$ or $8$ of the
$24$ units of $H_{1,2,2}$ are accounted for.
The complementary set of units appear when counting the number
of representations of certain multiples of an odd integer $m$
under parity restrictions.

\proclaim{Theorem 58 (Theorem on Complementary Representations)}
Let $m$ be an odd positive rational integer.
We consider representations by the quadratic form
$x^2 + y^2 + 2\, z^2 + 2\, w^2$.
\itemitem{(i)}
The number of representations of $4\, m$ by this quadratic form
with $x$ and $y$ even while $z$ and $w$ are odd is 
$4\, \sigma (m)$.
\itemitem{(ii)}
The number of representations of\/ $8 \, m$ by this 
quadratic form
where $x$ and $y$ are even while $z$ and $w$ are odd is 
$16 \, \sigma(m)$.
\itemitem{(iii)}
The number of representations of $4\, m$ by this quadratic form
where $x$ and $y$ are odd while $z$ and $w$ have different
parity is $16 \, \sigma (m)$.
\endproclaim

\demo{Proof}
Consider the equation in rational integers
$$
4 \, m \ = \ 
  x^2 + y^2 + 2\, z^2 + 2\, w^2 \,, \qquad
  x,\, y \hbox{ even,} \ \  
  z,\, w \hbox{ odd.}
\tag 7.5
$$
Dividing by $4$ we find
$w = (x/2)^2 + (y/2)^2 + 2\,(z/2)^2 + 2\, (w/2)^2$.
Set
$$
{x \over 2} \ = \ x^\prime \,, \quad
{y \over 2} \ = \ y^\prime \,, \quad
{z \over 2} \ = \ {1 \over 2} + z^\prime \,, \quad
{w \over 2} \ = \ {1 \over 2} + w^\prime \,.
\tag 7.6
$$
Then $x^\prime,\, y^\prime,\, z^\prime,\, w^\prime \in \Bbb Z$
and
$$
m \ = \ 
  N \left( x^\prime + y^\prime \,\qi 
           + \left( {1 \over 2} + z^\prime \right) \sqrt 2 \, \qj
           + \left( {1 \over 2} + w^\prime \right) \sqrt 2 \, \qk
   \right)
   .
\tag 7.7
$$
Note $\sqrt 2\, \qk = 2\, \bold v_4 - \bold v_1 - \bold v_2$,
$\sqrt 2\, \qj = 2\, \bold v_3 - \bold v_1 - \bold v_2$ and
$\bold v_3 + \bold v_4 = \bold v_1 + \bold v_2 
   + {{\sqrt 2} \over 2} (\qj + \qk)$.
The quaternion inside the norm symbol of (7.7) can be written
$$
\bold v_3 + \bold v_4 - \bold v_1 - \bold v_2
  + x^\prime \bold v_1 + y^\prime \bold v_2
  + z^\prime (2 \bold v_3 - \bold v_1 - \bold v_2 )
  + w^\prime (2 \bold v_4 - \bold v_1 - \bold v_2 )
  \,.
\tag 7.8
$$   
This is clearly an element of $H_{1,2,2}$ of norm $m$.
It has integer coefficients for $1$ and $\qi$ and half integer
coefficients for $\sqrt 2 \, \qj$ and $\sqrt 2 \, \qk$.

Any such element of $H_{1,2,2}$ of norm $m$ corresponds to a solution
of (7.5).
Let $\bold x$ be such an element, then $\bold x$ can be uniquely
factored as $\bold u \, \bold c$
for $\bold u$ a unit and $\bold c$ primary of norm $m$.
Arguing as in Theorem 57, $\bold x$ has the requisite coefficients
if and only if $\bold u$ has.
The number of units in $H_{1,2,2}$ with half integer coefficients
only for $\sqrt 2 \, \qj$ and $\sqrt 2 \, \qk$ is $4$, namely
$$
\pm \left\{
  {{\sqrt 2} \over 2} \qj \pm {{\sqrt 2} \over 2} \qk
  \right\}
\,.
\tag 7.9
$$
Thus the total number of solutions is $4\, \sigma (m)$.

Now consider the equation in rational integers
$$
8 \, m \ = \ 
  x^2 + y^2 + 2\, z^2 + 2\, w^2 \,, \qquad
  x,\, y \hbox{ even,} \ \  
  z,\, w \hbox{ odd.}
\tag 7.10
$$
The same substitution as in (7.6) yields the following variant
of (7.7):
$$
2\, m \ = \ 
  N \left( x^\prime + y^\prime \,\qi 
           + \left( {1 \over 2} + z^\prime \right) \sqrt 2 \, \qj
           + \left( {1 \over 2} + w^\prime \right) \sqrt 2 \, \qk
   \right)
   .
\tag 7.11
$$
As in the previous case, the term inside the norm symbol is an
element of $H_{1,2,2}$ with integer coefficients for $1$ and $\qi$
and half integer coefficients for $\sqrt 2 \, \qj$ and 
$\sqrt 2 \, \qk$.
Again, any such quaternion $\bold x$ is a solution to (7.10).
Such an $\bold x$ has a unique factorization as 
$(1 + \qi)$ times a unit $\bold u$ 
times a primary quaternion $\bold c$ of norm $m$.
As in Theorem 57, $\bold x$ has the requisite type of coefficients
if and only if $(1 + \qi) \, \bold u$ has.
Computation shows that the number of units where
$(1 + \qi) \, \bold u$ have half integer coefficients only for
$\sqrt 2 \, \qj$ and $\sqrt 2 \, \qk$ is $16$.
Thus the total number of solutions to (7.10) is $16 \, \sigma (m)$.

The remaining result breaks into two cases.
Consider the equation
$$
4 \, m \ = \ 
  x^2 + y^2 + 2\, z^2 + 2\, w^2 \,, \qquad
  x,\, y,\, z \hbox{ odd,} \ \  
  w \hbox{ even.}
\tag 7.12
$$
Dividing by $4$ and setting
$$
{x \over 2} \ = \ {1 \over 2} + x^\prime \,, \quad
{y \over 2} \ = \ {1 \over 2} + y^\prime \,, \quad
{z \over 2} \ = \ {1 \over 2} + z^\prime \,, \quad
{w \over 2} \ = \ w^\prime \,.
\tag 7.13
$$
we have $x^\prime,\, y^\prime,\, z^\prime,\, w^\prime \in \Bbb Z$
and
$$
m \ = \ 
  N \left( 
           {1 \over 2} + x^\prime 
	   + \left( {1 \over 2} + y^\prime \right) \,\qi 
           + \left( {1 \over 2} + z^\prime \right) \sqrt 2 \, \qj
           + w^\prime                              \sqrt 2 \, \qk
   \right)
   .
\tag 7.14
$$
The quaternion inside the norm symbol equals
$$
\bold v_3 
  + x^\prime \bold v_1 + y^\prime \bold v_2
  + z^\prime (2 \bold v_3 - \bold v_1 - \bold v_2 )
  + w^\prime (2 \bold v_4 - \bold v_1 - \bold v_2 )
  \,.
\tag 7.15
$$   
This is clearly an element of $H_{1,2,2}$ of norm $m$.
It has half integer coefficients for $1$, $\qi$ and $\sqrt 2 \, \qj$
but an integer coefficient for $\sqrt 2 \, \qk$.

The case of the equation
$$
4 \, m \ = \ 
  x^2 + y^2 + 2\, z^2 + 2\, w^2 \,, \qquad
  x,\, y,\, w \hbox{ odd,} \ \  
  z \hbox{ even, }
\tag 7.16
$$
is completely analogous, with solutions corresponding to 
those elements of norm $m$ with half integer coefficients for
$1$, $\qi$ and $\sqrt 2 \, \qk$  but an integer coefficient
for $\sqrt 2 \, \qj$.

Any $\bold x \in H_{1,2,2}$ can be uniquely written in the form
of a unit $\bold u$ times a primary quaternion $\bold c$ of norm
$m$.
Then as above, $\bold x$ has the coefficients of the proper
type if and only if $\bold u$ has.
Since the number of units with half integer coefficients
for $1$, $\qi$ and only one of $\sqrt 2 \, \qj$ or
$\sqrt 2 \, \qk$ is $16$, the total number of solutions
is $16 \, \sigma (m)$.
\locqed\enddemo

\bigskip
\head 8. Further Directions\endhead
It is very plausible that a similar study of the appropriate
norm Euclidean
quaternionic rings could lead to results for the
quadratic forms
$x^2 + 2y^2 + 3z^2 + 6w^2$ and
$x^2 + y^2 + 3z^2 + 3w^2$.
However, due to the significance of the prime $3$, the definition
of primary would have to be modified in an as yet undetermined
fashion.
Also the formula for the number of representations for the latter
quadratic form tends to imply that additional complications will
have to be surmounted.
This formula for the number of representations of $N$ is
$$
(-1)^{N-1} \, 4 \, \sum_{d \mid N} \, d \cdot \chi (d) \,,\qquad
  \chi(n) = \cases 
                      1   & \hbox{if } n \equiv 1,\, 5 \pmod 6     \cr
		      -1  & \hbox{if } n \equiv 2,\, 4 \pmod 6     \cr
		      0   & \hbox{if } n \equiv 0,\, 3 \pmod 6 \,. \cr
            \endcases 
\tag 8.1
$$
In particular the sum of divisors function does not appear in (8.1).
For a proof of the formula, see Fine [3].

\bigskip
\head 9. The Computation\endhead
The PUNIMAX version of MAXIMA was used on a LINUX partition 
of a Pentium 133 chip personal computer with 32 megabytes of RAM.
The operating system was Linux 2.0.35.
PUNIMAX was built using CLISP 1997-05-03 and GCC 2.7.2.3.
A short C program was written to numerically verify the Theorem
on Complementary Representations.

\bigskip
\head 10. Acknowledgments\endhead
The author would like to thank B.~Haible, the maintainer
of PUNIMAX, for permitting its free use for academic purposes 
[4].

\refstyle{C}

\enddocument